\input amstex
\input diagrams
\documentstyle{amsppt}

\NoRunningHeads \magnification=1000 \TagsOnRight \NoBlackBoxes
\hsize=5.3in \vsize=7.6in \hoffset= 0in \baselineskip=12pt
\topmatter
\title  Symmetric tensors and the geometry of subvarieties of $\Bbb P^N$.
\endtitle
\author Fedor Bogomolov*
 Bruno De Oliveira**
\endauthor
\thanks
* Partially supported by the NSF grant DMS-0100837.
** Partially supported by the NSF grant DMS-0306487
\endthanks
\affil
* Courant Institute for Mathematical Sciences
**University of Miami
\endaffil
\address
Fedor Bogomolov Courant Institute for Mathematical Sciences,
New York University\\
Bruno De Oliveira University of Miami
\endaddress
\email bogomolo{\@}CIMS.NYU.EDU bdeolive@math.miami.edu
\endemail
\endtopmatter
\document

\head {0. Introduction} \endhead

\

\

This paper following a geometric approach proves  new, and reproves old, vanishing and
nonvanishing results on the space of twisted symmetric differentials,
$H^0(X,S^m\Omega^1_X\otimes \Cal O_X(k))$ with $k\le m$, on subvarieties $X\subset \Bbb P^N$.
The case of $k=m$ is special and the nonvanishing results are related to the  space of quadrics
containing $X$ and lead to interesting geometrical objects associated to $X$, as for example the
variety of all tangent trisecant lines of $X$. The same techniques give results on the symmetric
differentials of subvarieties of abelian varieties. The paper ends with new results and examples about
the jump along smooth families of projective varieties $X_t$ of the symmetric plurigenera, $Q_m(X)=
\dim H^0(X,S^m\Omega^1_X)$, or of the $\alpha$-twisted symmetric plurigenera, $Q_{\alpha,m}(X)=
\dim H^0(X,S^m(\Omega^1_X\otimes \alpha K_X))$.

\

The paper is in part motivated by a previous result where we
proved that while  smooth hypersurfaces in $\Bbb P^3$ do not have symmetric differentials,
resolutions of nodal hypersurfaces  have them if the number of nodes is sufficiently large.
This is, in particular, interesting because smooth  and resolutions of nodal hypersurfaces
in $\Bbb P^3$of the same degree are deformation equivalent. So we have a case of the  jumping
of the symmetric plurigenera $Q_m(X)$ and a special one for that
matter, as we shall see below.
There is a previous example of this phenomenon in [Bo2-78], see section 2. Recall that this
contrasts with  the invariance of the plurigenera $P_m=\dim H^0(X,(\wedge^{n}
\Omega^1_X)^m)$, $n=\dim X$, [Si98]. The jumping in our example might bring back symmetric
differentials  to new approaches to the Kobayashi's conjecture which states that a general
hypersurace in $\Bbb P^3$ of degree $d\ge 5$ is hyperbolic (the known approaches use jet
differentials, i.e. higher order symmetric differentials, which exist on hypersurfaces but
are quite difficult to control). The authors motivated by this unexpected appearance of
the symmetric differentials  realized that there are many unanswered or forgotten questions
about them.

\

P. Bruckman showed in [Br71] that there are no symmetric differentials
on smooth hypersurfaces in $\Bbb P^N$ via an explicit constructive approach.
Later, F. Sakai with a cohomological approach using a vanishing theorem of
Kobayashi and Ochai showed that a complete intersection
$Y\subset \Bbb P^N$ with dimension $n>N/2$ has no symmetric differentials [Sa78].
In the early nineties M. Schneider [Sc92] using a similar approach, but  with
more general vanishing  theorems of le Potier, showed that any submanifold
$X\subset \Bbb P^N$ of dimension $n>N/2$ has no symmetric differentials of order $m$
even if twisted by $\Cal O_X(k)$, $H^0(X,S^m\Omega^1_X\otimes \Cal O_X(k))=0$, where $k<m$.

\

In this paper we use a distinct approach to obtain  vanishing and nonvanishing results
 on the space of twisted symmetric
differentials that are the natural
extension of  the results mentioned above.
Parts of this approach can be traced back to an announcement in the ICM of 1978
by the first author, [Bo78]. Our method has a geometric flavor involving
of the tangent map for $X$. The tangent map is given by $f:\Bbb P(\widetilde
{\Omega^1_X}(1)) \to \Bbb P^N$, where $f(\Bbb P(\widetilde {\Omega^1_X}(1))_x)
=T_xX$ and  $T_xX$ is the embedded projective tangent space to $X$ at $x$ in $\Bbb P^N$.
The first application of our approach  is to show that if the
tangent map for $X$  has a positive dimensional
general fiber, then $X$ has no symmetric differentials of order $m$ even if
twisted by $\Cal O_X(m)\otimes L$, where $L$ is any negative line bundle
on $X$, i.e. $H^0(X,S^m [\Omega^1_X(1)]\otimes L)=0 $ (this includes
Schneider's result).\par
The next step in the paper is to analyse the case of symmetric
powers of the  sheaf of twisted differentials , i.e.
$S^m[\Omega^1_X(1)]$. This would be  the case with  $L=\Cal O_X$ in
the result just mentioned above or equivalently the case $k=m$ not reached by
Schneider's methods. The loss of the negatitivity of $L$ makes the existence of
twisted symmetric differentials more delicate.  We use the pivotal  lemma 1.1 about
the sections of symmetric powers of quotients of trivial vector bundles. One of the
requirements to use the lemma is that the tangent map associated to $X\subset \Bbb P^N$
must be surjective and connected. A key result, theorem 1.3, is that connectedness of
the fibers of the tangent map is guaranteed if $dim X>2/3(N-1)$.\par

An, perhaps surprisingly, important feature of the existence of twisted symmetric differentials
is that it depends on the t-trisecant variety of $X$. The t-trisecant variety of
$X$ is the subvariety of the trisecant variety of $X$ consisting of the
union of the tangent lines to $X$ which are also trisecant lines. In fact, the
answer might also depend on of higher level t-trisecant varieties of $X$, see section
1.2.

We show in theorem C that if $X\subset \Bbb P^N$ has dimension
$n>2/3(N-1)$, then  $H^0(X,S^m [\Omega^1_X(1)])\neq 0 $ if and only if  all  higher
order t-trisecant varieties $X$ are not $\Bbb P^N$. This holds
in particular if $X$ is contained in a quadric.
When $X$ is  of codimension 1 or 2 and of dimension $n\ge 2$ $n\ge 3$, respectively,
the result is $H^0(X,S^m [\Omega^1_X(1)])\neq 0 $ if and
only $X$ is contained in a quadric (the codimension 2 case is the more challenging case).
Here an important characteristic is that one only needs to consider the t-trisecant variety.
Moreover, it is shown that in these cases the t-trisecant variety is the trisecant variety.
We also give a general criterion of when  the t-trisecant and the trisecant variety  of $X$
coincide. As an application of the circle ideas behind this criterion we give an alternative
proof of the Zak's theorem on the equality $Tan(X)=Sec(X)$ if  $Sec(X)$ does not have the
expected dimension.

\

\

In the last section we answer  a question of M. Paun asking if
$\dim H^0(X_t,S^m\Omega^1_{X_t}\otimes K_{X_t}))$ is locally
invariant in smooth families if $K_{X_t}>0$. We give a negative
answer based on the results on nonexistence twisted symmetric on
hypersurfacesthe plus the example, mentioned above, of the family
of smooth hypersurfaces in $\Bbb P^3$ specializing to the resolution
of a nodal surfaces with sufficiently many nodes. We then ask
whether there is a ratio $k/m$ for which the invariance of
$\dim H^0(X_t,S^m\Omega^1_{X_t} \otimes  (\wedge^{n} \Omega^1_X))^k)$
holds and  address the question of what is the lowest of such ratios.

\

\

\

\

 \head {1. Symmetric differentials on subvarieties
 of $\Bbb P^N$ and of Abelian varieties}\endhead

 \

 \

 \subhead  {1.1 Preliminaries}\endsubhead

 \

 \

Let $E$ be a vector bundle on $X$ and $\Bbb P(E)$ be the
projective bundle of hyperplanes of $E$. Recall the connection
between $S^mE$ and $\Cal O_{\Bbb P(E)}(m)$ which plays a
fundamental role in the study of symmetric powers of a vector
bundle. If $\pi: \Bbb P(E) \to X$ is usual projection map then the
following holds $\pi_*\Cal O_{\Bbb P(E)}(m)\cong S^mE$ and

$$H^0(X,S^mE)\cong H^0(\Bbb P(E),\Cal O_{\Bbb P(E)}(m)) \tag 1.1.1$$

The following case  persistently  appears in our arguments. Let
$E$ be a vector bundle on $X$  which is a quotient of
$\bigoplus^{N+1}L$ where $L$ is a line bundle on $X$:

$$q:\bigoplus^{N+1}L\to E \to 0 \tag 1.1.2$$

 Let $\Bbb
P(E)$ and $\Bbb P(\bigoplus^{N+1}L)$ be the  projective bundles of
hyperplanes of $E$ and $\bigoplus^{N+1}L$ respectively. The
surjection $q$ induces an inclusion and the isomorphism:

$$i_q:\Bbb P(E) \hookrightarrow \Bbb P(\bigoplus^{N+1}L)$$
$$ i_q^*\Cal O_{\Bbb
P(\bigoplus^{N+1}L)}(1)\cong \Cal O_{\Bbb P(E)}(1)$$

\noindent  Recall that there is a natural isomorphism $\phi: \Bbb
P(\bigoplus^{N+1} L)\to \Bbb P(\bigoplus^{N+1} \Cal O_X)$ for
which $\phi^*\Cal O_{\Bbb P(\bigoplus^{N+1}\Cal O)}(1)\cong \Cal
O_{\Bbb P(\bigoplus^{N+1}L)}(1)\otimes \pi^*L^{-1}$. The
projective bundle $\Bbb P(\bigoplus^{N+1} \Cal O_X)$ is the
product $X\times \Bbb P^N$, if
 $p_2$ denotes the projection onto the second factor, then $\Cal
O_{\Bbb P(\bigoplus^{N+1}\Cal O)}(1)\cong p_2^*\Cal O_{\Bbb
P^N}(1)$. Concluding, the surjection $q$ in (2) naturally induces
a map $f_q=p_2\circ \phi \circ i_q$ and the isomorphism:

$$f_q:\Bbb P(E) \to \Bbb P^N $$

$$f_q^*\Cal O_{\Bbb P^N}(1)\cong \Cal O_{\Bbb
P(E)}(1)\otimes \pi^*L^{-1}\tag 1.1.3$$

 Hence

$$H^0(X,S^mE)\cong H^0(\Bbb P(E),f_q^*\Cal O_{\Bbb P^N}(m)\otimes
\pi^*L^{\otimes m}) \tag 1.1.4$$

It follows from (1.1.4)  that the properties of the map $f_q: \Bbb
P(E) \to \Bbb P^N$ have an impact on the existence of sections of
the symmetric powers of $E$. The next result gives an example of
this phenomenon and will play a role in our study of existence of
symmetric differentials.

\proclaim {Lemma 1.1} Let $E$ be a vector bundle on a smooth
projective variety $X$. If $E$ is the quotient of a trivial vector bundle:

$$q:\bigoplus^{N+1}\Cal O_X\to E \to 0$$

\noindent and the induced map $f_q:\Bbb P(E) \to \Bbb P^N$ is
surjective with connected fibers, then $q$ induces the
isomorphism:

$$H^0(X,S^mE)= H^0(X,S^m(\bigoplus^{N+1}\Cal O_X))$$
\noindent  ($H^0(X,S^mE)=S^m[\Bbb Cs_0\oplus...\oplus \Bbb Cs_N]$ where $s_i=q(e_i)$,
$\bigoplus^{N+1}\Cal O_X=\bigoplus^{N}_{i=0}\Cal O_Xe_i$).
\endproclaim

\demo {Proof}  The isomorphism  $f_q^*\Cal O_{\Bbb P^N}(m)\cong \Cal O_{\Bbb P(E)}(m)$, (1.1.3), and
  $H^0(X,S^mE)\cong H^0(\Bbb P(E),\Cal O_{\Bbb P(E)}(m))$ give
  that:

 $$H^0(X,S^mE)\cong H^0(\Bbb P(E),f_q^*\Cal O_{\Bbb P^N}(m))$$

 The next step is to relate $ H^0(\Bbb P(E),f_q^*\Cal O_{\Bbb P^N}(m))$
 with   $H^0(\Bbb P^N,\Cal O_{\Bbb P^N}(m))$. If $f_q$ is
 surjective then $f_q^*:H^0(\Bbb P^N,\Cal O_{\Bbb P^N}(m))\to H^0(\Bbb P(E),f_q^*\Cal O_{\Bbb  P^N}(m))$
 is injective. If  the map $f_q$ also has connected fibers, then all
sections in $H^0(\Bbb P(E),f_q^*\Cal O_{\Bbb P^N}(m))$ descend to be sections in
$H^0(\Bbb P^N,\Cal O_{\Bbb P^N}(m))$, and the following holds:

$$H^0(\Bbb P(E),f_q^*\Cal O_{\Bbb P^N}(m))\cong H^0(\Bbb P^N,\Cal O_{\Bbb P^N}(m))$$

The result then follows from the  brake down of the map $f_q$, $f_q=p_2\circ i_q$,
plus $H^0(\Bbb P^N,\Cal O_{\Bbb P^N}(m)) \cong H^0(\Bbb P( \bigoplus^{N+1}\Cal O_X),p_2^*\Cal
O_{\Bbb P^N}(m))\cong H^0(\Bbb P( \bigoplus^{N+1}\Cal O_X),\Cal
O_{\Bbb P( \bigoplus^{N+1}\Cal O_X)}(m))$ and
$H^0(\Bbb P( \bigoplus^{N+1}\Cal O_X),\Cal O_{\Bbb P( \bigoplus^{N+1}\Cal O_X)}(m))
\cong H^0(X,S^m(\bigoplus^{N+1}\Cal O_X))$. \hfill \ \qed

\enddemo

\

\

\

\subhead  {1.2 Symmetric differentials on subvarieties of $\Bbb
P^N$}\endsubhead

\

\

The following is short collection of facts about the sheaf of differentials
that will help the reader understand our approach. The Euler sequence of $\Bbb P^N$ is:

$$0\to \Omega^1_{\Bbb P^N} \to \bigoplus^{N+1}\Cal O(-1) \to \Cal
O \to 0 \tag 1.2.1$$

\noindent The Euler sequence expresses the relation, induced by
the natural projection $p:\Bbb C^{N+1}\setminus \{ 0\} \to \Bbb
P^N$, between the differentials of $\Bbb C^{N+1}$ and $\Bbb P^N$.
A necessary condition for a differential $\omega$ of $\Bbb
C^{N+1}$ to come from a differential of $\Bbb P^N$ is that the
coefficients $h_0( z)$,...,$h_N( z)$ of $\omega=h_0(
z)dz_0+...+h_N( z)dz_N$ must be homogeneous  of degree -1. But the
last condition is not sufficient, the differentials $\omega$ on
$\Bbb C^{N+1}$ must be such that at any point $z\in \Bbb C^{N+1}$
their contraction with the vector $z_0\partial/\partial
z_0+...+z_N\partial/\partial z_N$, i.e with the direction of the
line from $z$ to the origin,  must be zero. To see this
algebraically,
 the sheaf $\bigoplus^{N+1}\Cal O(-1)$ in (1.2.1) is
$\bigoplus^{N+1}\Cal O(-1)=\Cal O(-1)dz_0+...+\Cal O(-1)dz_N$. The
map $q:\bigoplus^{N+1}\Cal O(-1) \to \Cal O \to 0$ is defined
sending $dz_i$ to $z_i$. So locally, let us say on $U_i=\{z_i\neq
0\}$, $\Omega^1_{U_i}$ the kernel of the map $q$ is spanned by the
sections induced by the differentials $\frac
{1}{z_i}dz_j-\frac{z_j}{z^2_i}dz_i$ on $p^{-1}(U_i)$.

\

The sheaf of differentials $\Omega^1_X$ is determined by (1.2.1)
restricted to $X$:

$$0\to \Omega^1_{\Bbb P^N}|_X \to \bigoplus^{N+1}\Cal O_X(-1) \to \Cal
O_X \to 0 \tag 1.2.2$$

\noindent and the conormal bundle exact sequence:

$$0\to N^* \to \Omega^1_{\Bbb P^N}|_X \to \Omega^1_X \to 0 \tag 1.2.3$$

The extension defined by (1.2.2) (which corresponds to a cocycle
$\alpha \in H^1(X,\Omega^1_{\Bbb P^N}|_X)$) induces via the
surjection in (1.2.3) the extension:

$$0\to \Omega^1_X  \to \widetilde {\Omega^1_X} \to \Cal
O_X \to 0 \tag 1.2.4$$

\noindent The geometric description of the sheaf
$\widetilde {\Omega^1_X}$ is that it is the sheaf on $X$ associated
to the sheaf of 1-forms on the affine cone $\hat X\subset \Bbb C^{N+1}$.
The above exact sequences after twisted by $\Cal O_X(1)$ fit
in the commutative diagram:

\diagram &&&0&\rTo&\Omega^1_{\Bbb
P^N}|_X(1)&\rTo&\bigoplus^{N+1}\Cal O_X&\rTo&\Cal O_X(1)&\rTo&0\\
&&&&&\dOnto&&\dOnto^{q}&&\dTo^{\simeq}&&&
\text{(1.2.5)}\\&&&0&\rTo&\Omega^1_X(1)&\rTo&\widetilde {\Omega^1_X}(1)
&\rTo&\Cal O_X(1)&\rTo&0
\enddiagram

\

The middle vertical
surjection of diagram (1.2.5) can be represented more explicitly by:

$$q:\bigoplus^{N}_{i=0}\Cal O_Xdz_i\to \widetilde {\Omega^1_X}(1)
 \tag 1.2.6$$

\noindent The induced map $f:\Bbb P(\widetilde {\Omega^1_X}(1))
\to \Bbb P^N$ is such that for each $x\in X$:

$$f(\Bbb P(\widetilde {\Omega^1_X}(1))_x)=T_xX \tag 1.2.7$$

\noindent where $T_xX$ is the embedded projective tangent space
to $X$ at $x$ inside $\Bbb P^N$. For the obvious reasons $f$ will
be called the tangent map for $X$. The tangent map $f$ induces a map
from $X$ to $G(n,N)$ which is exactly the Gauss map for $X$,
$\gamma_X:X\to G(n,N)$.

\

\proclaim {Theorem A} Let $X$ be a smooth projective subvariety of
$\Bbb P^N$. If the general fiber of the tangent map for $X$, $f:\Bbb P(\widetilde {\Omega^1_X}(1))
\to \Bbb P^N$, is positive dimensional, then $\forall m\ge 0$:

$$H^0(X,S^m[\Omega^1_X(1)]\otimes L)=0$$

if $L$ is a negative line bundle on $X$.
\endproclaim

\demo {Proof}  It is sufficient to show that $ H^0(X, S^m[\widetilde
 {\Omega^1_X}(1)]\otimes L)=0$, since there is the inclusion $S^m[\Omega^1_X(1)]
\hookrightarrow  S^m[\widetilde {\Omega^1_X}(1)]$,
induced from (1.2.5). \par

The projective bundle $\Bbb P(\widetilde {\Omega^1_X}(1))$ comes
with two maps. The tangent map for $X$, $f:\Bbb P(\widetilde
{\Omega^1_X}(1)) \to \Bbb P^N$, and the projection onto $X$,
$\pi:\Bbb P(\widetilde {\Omega^1_X}(1)) \to X$. One also has
the natural isomorphisms $\Cal O_{\Bbb P(\widetilde {\Omega^1_X}
(1))}(m)=f^*\Cal O_{\Bbb P^N}(m)$ and $\pi_*(\Cal O_{
\Bbb P(\widetilde{\Omega^1_X}(1))}(m)\otimes\pi^*L)
\cong S^m[\widetilde {\Omega^1_X}(1)]\otimes L$. These
isomorphisms give:

$$ H^0(X, S^m[\widetilde {\Omega^1_X}(1)]\otimes L)\cong
H^0(\Bbb P(\widetilde {\Omega^1_X}(1)), f^* \Cal O_{\Bbb P^N}(m)
\otimes \pi^*L)$$

\noindent The vanishing of the last group follows from the negativity
of the line bundle $f^* \Cal O_{\Bbb P^N}(m)\otimes \pi^*L$ along
each fiber of the map $f$. More precisely, $f^*\Cal O_{\Bbb P^N}(m)$
is trivial on the fibers and $\pi^*L$ is negative on the fibers since
$L$ is negative on $X$ the map $\pi$ is injective on each fiber of $f$.
\par
We need the fibers of the map $f$ to be positive dimensional. Since is
only in this case that the negativity of the line bundle
$f^* \Cal O_{\Bbb P^N}(m)\otimes \pi^*L$, $l<0$, makes sense.
This negativity implies that all sections of $H^0(\Bbb P(\widetilde
{\Omega^1_X}(1)), f^* \Cal O_{\Bbb P^N}(m)\otimes \pi^*L)$ vanish
along all fibers of $f$ and hence vanish on all $\Bbb P(\widetilde
{\Omega^1_X}(1))$, which completes the proof.
\hfill \ \qed
\enddemo

\

As an important  case of theorem A one has another proof to the result
first proved by Schneider [Sc92].

\

\proclaim {Corollary 1.2} Let $X$ be a smooth projective subvariety of
$\Bbb P^N$ whose dimension $n>N/2$. Then:

$$H^0(X,S^m\Omega^1_X\otimes \Cal O(k))=0$$

if $k<m$.
\endproclaim

\demo {Proof} The dimensional hypothesis $n>N/2$ guarantee that all fibers of
the tangent map $f$ for $X$ are positive dimensional. The condition $k<m$ gives that
$S^m\Omega^1_X\otimes \Cal O(k)=S^m[\Omega^1_X(1)]\otimes \Cal O_X(l)$, with $l<0$.
The theorem then follows from theorem A for the negative line bundle $L=\Cal O_X(l)$,
$l<0$.
\hfill \ \qed
\enddemo

\

What happens in the key case $k=m$? The results just mentioned use
the negativity $f^* \Cal O_{\Bbb P^N}(m) \otimes \pi^*L$,  along the fibers
of the map $f$, which no longer holds if $k=m$. Indeed, one has $H^0(X,
S^m\widetilde {\Omega^1_X}\otimes \Cal O_X(m))=H^0(\Bbb P(\widetilde
{\Omega^1_X}(1)), f^* \Cal O_{\Bbb P^N}(m))$ which is no longer trivial.
The  analysis of the nonexistence of twisted symmetric differentials
$\omega\in H^0(X,S^m[\Omega^1_X(1)])$ on $X$ is more delicate.  One has to describe
the sections  $H^0(X, S^m\widetilde {\Omega^1_X}\otimes \Cal O_X(m))$ and
characterize which ones are in $H^0(X, S^m {\Omega^1_X}\otimes \Cal O_X(m))$.
The answers will depend on  geometric properties involving the variety of
tangent lines to the subvariety $X$.

\

To describe the twisted symmetric extended differentials in
$H^0(X, S^m[\widetilde {\Omega^1_X}(1)])$ one needs to use the properties
of the tangent map for $X\subset \Bbb P^N$. The lemma 1.1 gives a good
description of $H^0(X, S^m[\widetilde{\Omega^1_X}(1)])$ if the map
tangent map $f$ is a connected surjection. The next lemma shows that
this is the case when $\dim X>2/3(N-1)$.

\

This paragraph  about the tautological $\Bbb P^1$-bundle over the
grassmanian is used in the lemma below. Let $q:I \to G(1,\Bbb
P^n)$ be the tautological $\Bbb P^1$-bundle over the grassmanian
and  $p:I \to \Bbb P^n$ the natural map. For any point $x\in \Bbb
P^n$ there is a $\Bbb P^{n-1}\subset G(1,\Bbb P^n)$ consisting of
all lines passing through $x$. The restriction (or the pullback)
of the tautological $\Bbb P^1$-bundle to $\Bbb P^{n-1}$ is $q:\Bbb
P(\Cal O(1) \oplus \Cal O) \to \Bbb P^{n-1}$.

\

\proclaim {Theorem 1.3} Let $X$ be a $n$-dimensional submanifold
of $\Bbb P^N$ with $n>2/3(N-1)$ then the natural map $f:\Bbb
P(\widetilde {\Omega^1_X}(1))\to \Bbb P^N$ associated with the
Gauss map is a surjective and connected morphism.
\endproclaim

\demo {Proof} The description of the map $f$ in (1.2.7) implies
that:

 $$f(\Bbb P(\widetilde {\Omega^1_X}(1)))=Tan(X)$$
\noindent where $Tan(X)$ is the tangent variety of $X$, in other
words $Tan(X)=\bigcup_{x\in X}T_xX\subset \Bbb P^N$. Denote by
$Sec(X)\subset \Bbb P^N$ the secant variety of $X$. Zak's results
[Za81] about tangencies state that one of the following must hold:
i) $\dim Tan(X)=2n$ and $\dim Sec(X)=2n+1$; ii) $Tan(X)=Sec(X)$.
It follows immediately that if $\dim X\ge N/2$ then
$Tan(X)=Sec(X)$.\par It is also a result of Zak, coming from
applying the results on tangencies plus the Terracini's lemma on
the tangent spaces of secant varieties, that $Sec(X)\neq \Bbb P^N$
can hold only if $n\le 2/3(N-2)$. Hence surjectivity of $f$ is
guaranteed if $n>2/3(N-2)$, which is the case.

\

It remains to show the connectedness of the fibers. Denote the
fibers of $f$ by $Y_x=f^{-1}(x)$ for $x\in \Bbb P^N$ and $\pi:\Bbb
P(\widetilde {\Omega^1_X}(1))\to X$ be the projection map. The
injectivity of $f$ restricted to the fibers of $\pi$ implies that
$Y_x$ is connected if and only if $\pi(Y_x)$ is connected. The
subvariety $R_x=\pi(Y_x)$ is the locus of $X$ consisting of all
the points in $X$ having a tangent line passing through $x$. The Stein factorization
implies that $f$ is connected if its general fiber is connected, i.e.
if for the general $x\in \Bbb P^N$ the locus $R_x$ is connected.\par

In the following arguments we always assume that $x\in \Bbb P^N$
is general.  The first observation to make is that $R_x\subset
Z_x$, where $Z_x$ is the double locus of the projection $p_x: X\to
\Bbb P^{N-1}$ (i.e. the locus of points in $X$ belonging to  lines
passing through $x$ and meeting $X$ at least twice). By
dimensional arguments one has that $R_x$ is a Weyl divisor of
$Z_x$. A key element in our argument is the result of [RaLo03]
stating that the double locus $Z_x$ is irreducible if
$n>2/3(N-1)$.
\par

Let $S\subset \Bbb P^{N-1}$ be the image of $Z_x$ by the
projection $p_x$. The irreducible variety $S$ can be seen as a
subvariety of the $\Bbb P^{n-1}\subset G(1,\Bbb P^N)$ of lines
passing through $x$. We can pullback the tautological $\Bbb
P^1$-bundle on $G(1,\Bbb P^N)$ to $S$ and obtain $q_S:\Bbb P(\Cal
O(1)|_S \oplus \Cal O) \to S$. The natural map $p:I \to \Bbb P^N$,
see the paragraph before the lemma, induces a map $p:\Bbb P(\Cal
O(1)|_S \oplus \Cal O) \to \Bbb P^N$, whose image is the cone with
vertex p and base $S$. The map $p$ is a biregular morphism of the
complement of $p^{-1}(x)$ onto the cone without the vertex.

\

The $\Bbb P^1$-bundle $q_S:\Bbb P(\Cal O(1)|_S \oplus \Cal O) \to
S$ comes with two natural sections (one for each surjection onto
the lines bundles $\Cal O$ and $\Cal O(1)$). Geometrically these
two sections come from the pre-image of $x$ and the pre-image of
$S$ via the map $p$. The subvariety $M= p^{-1}(Z_x)$ is biregular
to $Z_x$ and is a divisor in the total space  of the line bundle
$\Cal O_S(1)$, $\Bbb P(\Cal O_S(1)\oplus \Cal O)\setminus \Bbb
P(\Cal O)$. The points in $p^{-1}(R_x)$ are the points $y \in M$
for which the fibers of $q_S$ meet $M$ at $y$ with multiplicity $\ge 2$.
The generality of $x$ implies by the classical trisecant lemma
that the general fiber of $q_S$ meets $M$ only twice counting with
multiplicity. This makes the projection $q_S|_M:M \to S$ a generically 2 to 1 map.

\par

Consider the pullback $L=q_S|_M^*\Cal O_S(1)$ which is an ample
line bundle on $M$.  The line bundle $L$ comes naturally with a
nontrivial section denote the corresponding divisor of the total
space of $L$, $Tot(L)$, by $D_1$. Denote the natural map between
the total spaces of $L$ and $\Cal O_S(1)$ by $g:Tot(L) \to
Tot(\Cal O_S(1))$. The divisorial component of $g^{-1}(M)$ is
decomposed in two irreducible components $D_1$ and $D_2$. Let
$h:Tot(L) \to M$ be the natural projection, then $h(D_1\cap
D_2)\subset R_x$. If $D_2$ is also a section of $L$, then
$h(D_1\cap D_2)$ is connected since it is the zero locus $(s)_0$
of a section $s$ of the ample line bundle $L$. The result would follow
since any other possible component of $R_x$ has to meet $(s)_0$.
If $D_2$ is not a section the result stills follows from the same argument
after base change (pulling back $L$ to $D_2$ using $h$).

\hfill \ \qed
\enddemo

\

In conjunction  with lemma 1.1 one obtains the following description
of the space of  twisted extended symmetric differentials on $X$:

\

\proclaim {Corollary 1.4} Let $X$ be a $n$-dimensional submanifold of
$\Bbb P^N$ with $n>2/3(N-1)$ then:

$$H^0(X, S^m[\widetilde {\Omega^1_X}(1)])=S^m[\Bbb Cdz_0\oplus ...\oplus \Bbb Cdz_N].$$

\endproclaim

\

\noindent The characterization of the space of twisted symmetric
differentials on $X$, within the dimensional range
$\dim X>2/3(N-1)$, is given by the following proposition:

\

\proclaim {Proposition 1.5} Let $X$ be a $n$-dimensional submanifold
of $\Bbb P^N$ with $n>2/3(N-1)$ then:

$$H^0(X, S^m[ {\Omega^1_X}(1)])=\{\Omega \in
S^m[\Bbb Cdz_0\oplus ...\oplus \Bbb Cdz_N]|
 \text { } Z(\Omega)\cap T_xX\text { is a cone with vertex at } x,
  \text { } \forall x \in X\}$$
\endproclaim

\demo {Proof} The inclusion $H^0(X, S^m[ {\Omega^1_X}(1)])
\subset H^0(X, S^m[\widetilde {\Omega^1_X}(1)])$ and corollary  1.4
imply that that all the sections of $H^0(X, S^m[ {\Omega^1_X}(1)])$
are induced from the symmetric $m$-differentials $S^m[\Bbb Cdz_0
+...+\Bbb Cdz_N]$ on $\Bbb C^{N+1}$.  \par

Let $\hat X \subset \Bbb C^{N+1}$ be the affine cone over $X\subset
\Bbb P^N$, $T\hat X$ the sheaf on $X$ associated with the tangent
bundle of $\hat X$ and  $T_xX\subset \Bbb P^N$ the embedded tangent
space to $X$ at $x$. Consider the rational map $p:\Bbb P(\widetilde
{\Omega^1_X}(1)) \dashrightarrow \Bbb P({\Omega^1_X}(1))$,
which is fiberwise geometrically described by the projections from
the point $x \in T_xX$ $p_x:T_xX=\Bbb P_l(T\hat X)_x\dashrightarrow
\Bbb P_l(TX)_x$, ($\Bbb P_l(E)$ is the projective bundle of lines in
 the vector bundle $E$, $\Bbb P_l(E)=\Bbb P(E^*)$).
The map $p$ gives an explicit inclusion $H^0(X,S^m[{\Omega^1_X}(1)])
=p^*H^0(\Bbb P({\Omega^1_X}(1)), \Cal O_{\Bbb P( {\Omega^1_X}(1))}(m))
\subset H^0(\Bbb P(\widetilde {\Omega^1_X}(1)), \Cal O_{\Bbb P(\widetilde
{\Omega^1_X}(1))}(m))=S^m[\Bbb Cdz_0+...+\Bbb Cdz_N]$.\par
Recall that if $0\to V\to \widetilde V \to \Bbb C \to 0$ is a sequence
of vector spaces, then one gets a projection from $[V]\in
\Bbb P(\widetilde V)$, $p:\Bbb P(\widetilde V)\dashrightarrow \Bbb P(V)$.
The sections in $H^0(\Bbb P(\widetilde V),\Cal O_{\Bbb P(\widetilde V)}(m))$
which are in $p^*H^0(\Bbb P(V),\Cal O_{\Bbb P(V)}(m))$
are the ones corresponding to homogeneous polynomials whose zero locus
is a cone with vertex at $[V]$.\par
An element $\Omega \in S^m[\Bbb Cdz_0\oplus ...\oplus \Bbb Cdz_N]$
corresponds in a natural way to an homogeneous polynomial in $\Bbb P^N$
which we still denote by $\Omega$. From the last two paragraphs it
follows that $\Omega$ induces an element in $\omega \in H^0(X, S^m
[ {\Omega^1_X}(1)])$ if and only if $\forall x \in X$ the zero locus
$Z(\Omega)\cap {T_xX}$ is a cone with vertex  $x$.\hfill \ \qed
\enddemo

\

\

We proceed to extract from proposition 1.5  the geometric conditions
required for the existence of twisted symmetric differentials on
smooth subvarieties $X\subset \Bbb P^N$. First, we  need to introduce
some objects and notation.

\

Let $X$ be an irreducible subvariety and $Y$ be any subvariety of
$\Bbb P^N$. Consider the incidence relation:

$$\Cal C_XY:=\overline{\{(x,z)\in X_{sm}\times \Bbb P^N \text { }|
 \text { } z\in \overline {xy},\text { } y\neq x \text { and } y
 \in Y\cap T_xX\}  }\subset X\times \Bbb P^N$$

\noindent The variety $\Cal C_XY$ comes with two projections. Denote by
$C_XY:=p_2(\Cal C_XY)$. Equivalently, $\forall x\in X_{sm}$ denote by
$C_xY\subset T_xX$ the cone with vertex at $x$ consisting of the closure
of the union of all chords joining $x$ to $y\neq x$ with $y\in Y\cap T_xX$,
where $T_xX$ is the projective embedded tangent space to $X$ at $x$. Then
$C_XY=\overline {\bigcup_{x\in X_{sm}}C_xY}\subset \Bbb P^N$.

\

\proclaim {Definition 1.6} Let $X$ be an irreducible subvariety of
$\Bbb P^N$. The variety $C_XX\subset \Bbb P^N$ will be called the
t-trisecant variety of $X$.
\endproclaim

\

For hypersurfaces $H$, including singular but always reduced, one
has the following useful result.

\

\proclaim {Proposition 1.7} Let $H\subset \Bbb P^N$, $N> 2$, be an irreducible
nondegenerate  hypersurface. Then  $C_HH=H$ if $H$ is a
quadric and  $C_HH=\Bbb P^N$ otherwise.
\endproclaim

\demo {Proof}   Let $H$ be a quadric and $x\in H$
a smooth point. The lines $l\subset C_xH$ passing through $x$ must
touch $H$ at least 3 times (counting with multiplicity) hence
$l\subset H$. This implies that $C_xH\subset H$, for all $x\in H$,
therefore $C_HH=H$.\par
Let $H$ be of degree greater than 2. The result follows to the trivial case
of curves in $\Bbb P^2$. For a general $2$-plane $L$ in $\Bbb P^N$ the intersection
$H\cap L=D$ is an irreducible and reduced curve of the same degree as $H$.
The irreducible and reduced curve $D=H\cap L$ of degree $\ge 3$ in $L=\Bbb P^2$
satisfies $C_DD=L$ (well known but see remark below). The result follows since
the following inclusion holds $C_{(H\cap L)}(H\cap L)\subset (C_HH)\cap L$ and
hence $C_HH$ contains the general 2-plane.

\hfill \ \qed
\enddemo

\

\

\proclaim {Remark} Let $x\in H$ be a general point and
$\Bbb C^N\subset P^N$ be an affine chart containing   $x$,
where w.l.o.g.  $x=0$. Let the hypersurface $H\cap \Bbb C^N$
be  given by  ${f=0}$. The quadractic part $f_2$ of the Taylor
expansion of $f$ at $x$ can not be trivial on $T_xH$, otherwise
there would be an open subset of $H$ on which the second
fundamental form of $H$ is trivial which would force $H$ to be
an hyperplane. Denote by $Q_x$ the quadric defined by $f_2|_{T_xH}$.
In the case $H\subset \Bbb P^2$ of degree $d\ge 3$, then $Q_x$ not
being trivial implies that tangent line $l=T_xH$ is such that
$(l\cap H)_x=2$ and thus $l$ must meet $H$ away from $x$. Hence
$C_HH=Tan(H)$.
\endproclaim

\

\

The first appearance of the t-trisecant variety
is on the  result about the existence of twisted symmetric differentials
on smooth hypersurfaces in $\Bbb P^N$.

\
\proclaim {Theorem B} Let $X$ be a smooth projective hypersurface
in $\Bbb P^N$. Then:

$$H^0(X,S^m[\Omega^1_X(1)])=0$$

\noindent if and only if the $C_XX=\Bbb P^N$, i.e. X is a not quadric.
\endproclaim
\demo {Proof} Let $\Omega\in S^m[\Bbb Cdz_0\oplus ...\oplus \Bbb Cdz_N]$ be
such that it induces a nontrivial $\omega \in H^0(X, S^m[ {\Omega^1_X}(1)])$.
As before denote  by $\Omega$ also the corresponding homogeneous polynomial
of degree $m$ in $z_0,...,z_N$. Proposition 1.5 says that the zero locus
$Z(\Omega)$ is such that $\forall x \in X$ $Z(\Omega)\cap T_xX$ is a cone with
vertex $x$. This clearly implies that $X\subset Z(\Omega)$. Moreover since
for $ x \in X$ $Z(\Omega)\cap T_xX$ is a cone with vertex $x$ containing
$X\cap T_xX$, then $C_xX\subset Z(\Omega)\cap T_xX$. Hence:

$$C_XX\subset Z(\Omega) \tag 1.2.8$$

\noindent Therefore if $C_XX=\Bbb P^N$ then $\Omega=0$ and the ($\Leftarrow$)
part of the theorem is proved.

\

If $C_XX\neq \Bbb P^N$ then $X\cap T_xX$ is a cone with vertex $x$ for the
general $x\in X$. To see this notice that if $X\cap T_xX$ is not a cone with
vertex $x$ then by dimensional arguments $C_xX=T_xX$. Since $Tan(X)=\Bbb P^N$
for hypersurfaces, $C_XX\neq\Bbb P^N$ implies that  $C_xX\neq T_xX$ for the
general $x\in X$. The condition that $X\cap T_xX$ is a cone with vertex $x$ is
a closed condition on $x\in X$, hence $C_XX=X$ if $C_XX\neq \Bbb P^N$.\par
If $C_XX=X$ then any differential $\Omega\in S^m[\Bbb Cdz_0\oplus ...\oplus
\Bbb Cdz_N]$ whose zero locus is a multiple of the hypersurface $X$ will induce
a nontrivial  element $\omega \in H^0(X, S^m[ {\Omega^1_X}(1)])$.
Proposition 1.7 states that if $C_XX=X$ if and only if $X$ is a quadric and
finishes the proof.
\hfill \ \qed
\enddemo

\

\proclaim {Remark} Standard arguments can easily show that quadrics have
twisted symmetric differentials. For example, let $X\subset \Bbb P^3$ be
a nonsingular quadric. Then the  surface $X$ is $\Bbb P^1\times \Bbb P^1$,
$\Omega^1_X=\Cal O_{\Bbb P^1\times \Bbb P^1}(-2,0)\oplus \Cal O_{\Bbb P^1
\times \Bbb P^1}(0,-2)$ and $\Cal O_X(1)=\Cal O_{\Bbb P^1\times \Bbb P^1}(1,1)$.
Hence $S^2[\Omega^1_X(1)]=\Cal O_{\Bbb P^1\times \Bbb P^1}(-2,2)\oplus
\Cal O_{\Bbb P^1\times \Bbb P^1}\oplus \Cal O_{\Bbb P^1\times \Bbb P^1}(2,-2)$
which implies  $H^0(X,S^2[\Omega^1_X(1)])=\Bbb C$ as expected from
proposition 1.5 or theorem B.
\endproclaim

\

\

We proceed to analyse the higher codimension case for $X\subset
\Bbb P^N$. In the hypersurface case the knowledge of the t-trisecant variety
$X$ was sufficient to obtain the complete answer in
theorem B. But in higher codimension,  one should also consider iterations of
the construction $C_XY$. Define  $C_X^2Y=C_X(C_XY)$ (note $C_X^2Y
\neq C_{_{(C_{_X}Y)}}(C_XY)$) and proceed inductively to obtain
$C^k_XY$.

\

\proclaim {Theorem C} Let $X$ be a non degenerated smooth
projective subvariety of $\Bbb P^N$ of dimension $n> 2/3(N-1)$. If
$C_X^{k}X=\Bbb P^N$ for some $k$, Then:

$$H^0(X,S^m[\Omega^1_X(1)])=0$$

\endproclaim

\demo {Proof} It follows from the proposition 1.5 that in  the dimensional range $n > 2/3(N-1)$  the differentials
$\omega \in H^0(X, S^m[ {\Omega^1_X}(1)])$ are induced from symmetric
$m$-differentials $\Omega \in S^m[\Bbb Cdz_0+...+\Bbb Cdz_N]$ on
$\Bbb C^{N+1}$. Moreover, Proposition 1.5 also says that the zero locus
$Z(\Omega)$  must be such that  $Z(\Omega)\cap T_xX$ is a cone with
vertex $x$, $\forall x \in X$. In the proof of theorem B, it was shown
that this implies  that $C_XX\subset Z(\Omega)$. \par Following the same
reasoning, since $Z(\Omega)\cap T_xX$ is a cone with vertex $x$ and
$C_XX\subset Z(\Omega)$ then $C_X^2X\subset Z(\Omega)$. Repeating the
argument one gets $C_X^lX\subset Z(\Omega)$ for all $l\ge 1$. If $C_X^kX=
\Bbb P^N$ for some $k$, then clearly every symmetric differential
$\Omega\in S^m[\Bbb Cdz_0+...+\Bbb Cdz_N]$ inducing $\omega \in H^0
(X, S^m[ {\Omega^1_X}(1)])$ must be trivial.
\hfill \ \qed
\enddemo

\

\

It follows from theorem B and C that it is important to characterize
the subvarieties $X\subset \Bbb P^N$ with $\dim X>2/3(N-1)$ with $C^k_XX
\neq \Bbb P^N$ for any $k$. They must be  special and, as in the hypersurface
case, subjectable to a description. A general result showing again
the role of quadrics is:

\

\proclaim {Proposition 1.8} Let $X$ be a subvariety of $\Bbb P^N$ such
that $X\subset Q_1\cap ...\cap Q_l$, where $Q_1$,..., $Q_l$ are quadrics.
Then $C^k_XX\subset Q_1\cap ...\cap Q_l$ for all $k\ge 1$.
\endproclaim

\demo {Proof} There is the  inclusion of the t-trisecant varieties
$C_XX\subset C_{Q_i}Q_i$ for all quadrics $Q_i$ $i=1,...,l$.  The
equality $C_{Q_i}Q_i=Q_i$ proved in proposition 1.6 gives $C_XX
\subset Q_1\cap ...\cap Q_l$. In an  equal fashion one sees that
$C_X^2X=C_X(C_XX)\subset C_{Q_i}Q_i$ for all $i=1,...,l$.
Induction then gives the result.\hfill \ \qed
\enddemo

\

\proclaim {Remark} One should investigate  the conditions on $X$ for
the validity of the assertion that  $C^k_XX$ is the intersection of all
quadrics containing $X$ for $k$ sufficiently large.
\endproclaim

\

The case where $X$ is of
codimension 2, $X^n\subset \Bbb P^{n+2}$, also has a complete answer, see
proposition 1.12 below. The answer it will follow from establishing  that $C_XX$ is the
trisecant of $X$, $Tr(X)$, if $n\ge 3$, then  one can use the results  on
trisecant varieties of varieties of codimension 2 of Ziv Ran [Ra83],
$n\ge 4$, and Kwak [Kw02] for the threefold case.

\

Let $X\subset \Bbb P^N$ be a subvariety and $l\subset \Bbb P^N$  a line meeting
$X$ at $k$ points, $x_i$  $i=1,...,k$, the line $l$ is said to is of type
$(n_1,...,n_k)$ if $n_i=length_{x_i}(X\cap l)$. A line $l$ is a trisecant line if
$\sum n_i\ge 3$ and a t-trisecant line if additionally one of the $n_i\ge 2$ ($Tr(X)$
is the union of all trisecant lines and $C_XX$ is the union of all t-trisecant lines).

\proclaim {Lemma 1.9} Let $X$ be a subvariety of $\Bbb P^N$ and $\pi:L \to T$
a 1-dimensional family of lines in $\Bbb P^N$ all passing through a fixed $z\notin X$ and
whose union is not a line. If the general line meets $X$ at least twice, then one
of the lines must meet $X$ with multiplicity at least 2 at some point.
\endproclaim

\demo {Proof} Let $H$ be an hyperplane not containing $z$ and $f:T \to H$
be the map which sends $t$ to $L_t\cap H$. Denote by $C$ the image of map
$f$. Let $C(z,C)$ the cone over $C$ with vertex $z$. Let $D$ be the curve
which consists of the divisorial component of $X\cap C(z,C)$. The possibly
nonreduced curve $D$ is such that any line $l_{c}$ joining $z$ to $c\in C$ meets
$D$ at least twice (counting with multiplicity). We can  assume that
the lines $l_c$ meet  $X$ with at most multiplicity at any $x\in X$,
otherwise the result follows. Hence the the curve $D\subset C(z,C)$ is reduced
and clearly does not pass through $z$.\par

Resolve the cone $C(z,C)$ by normalizing $C$, $\bar C$, and
blowing up the singularity at the vertex. The resulting surface
$Y$ is a ruled surface over $\bar C$, which comes with two
maps $\sigma: Y\to C(z,C)$ and $f:Y \to \bar C$. Let $\bar D$
be the pre-image of $D$ by $\sigma$. If $\bar D$ meets any of the
fibers of $f:Y \to \bar C$ with multiplicity $\ge 2$ then as before
we are done. Hence $\bar D$ is smooth moreover it must be a
multi-section. This is impossible since by base change we would
obtain a ruled surface which would have at least two disjoint
positive sections not intersecting (the unique
negative section lies over the the pre-image of $p$).
\hfill \ \qed
\enddemo

We proceed by giving an alternative proof of Zak's theorem on the equality of
the secant and tangent variety for smooth subvarieties $X$ whose secant variety
does not have the expected dimension.

\

\proclaim {Corollary 1.10} (Zak's Theorem) Let $X$ be a smooth subvariety
of $\Bbb P^N$. If $\dim Sec(X)<2n+1$ then $Tan(X)=Sec(X)$.
\endproclaim

\demo{Proof} Assume $Sec(X)\neq X$ since if the equality holds then clearly
$Tan(X)=Sec(X)$. Let $z$ be a point of $Sec (X) \setminus X$.
Since $Sec(X)$ has less than the expected dimension there is  a positive dimensional
family $\pi:L \to T$ of secant lines passing through $z$. Apply lemma 1.10 to
a 1-dimensional subfamily of $\pi:L \to T$ and obtain that one of this lines $L_{t_0}$
must meet $X$ with multiplicity at least 2 at some point $x\in X$ hence $L_{t_0}$ is tangent to $X$ at $x$
and $z\in Tan(X)$.
\hfill \ \qed
\enddemo

\

Finally we use the lemma to describe an important case when the
trisecant variety is equal to the t-trisecant variety.

\

\proclaim {Corollary 1.11} Let $X$ be a smooth subvariety
of $\Bbb P^N$.  If  the family of  trisecant lines of  $X$ through
a general point of $Tr(X)$ is at least 1-dimensional, then $Tr X=C_XX$.
\endproclaim

\demo {Proof} The same argument after replacing $Sec(X)$ by $Tr(X)$ and
$Tan(X)$ by $C_XX$.
\hfill \ \qed
\enddemo

\

\

\

\proclaim {Proposition 1.12} Let $X$ be a smooth subvariety of codimension 2
in $\Bbb P^{n+2}$. If $n\ge 3$ then:

\

 1) $C_XX=Tr(X)$.

 \

 2)   $C_XX=\Bbb P^{n+2}$ or
$C_XX$ is the intersection of the quadrics containing $X$.
\endproclaim

\demo {Proof} First we establish 1). Let $z$ be a general point
of the trisecant variety $Tr(X)$. Let $l$ be a trisecant line passing
through $z$, assume it is not also t-trisecant since otherwise
there is nothing to prove. Consider the projection $p_z:X
\to \Bbb P^{n+1}$ from the point $z$ to an hyperplane $\Bbb P^{n+1}
\subset \Bbb P^{n+2}$. Denote 3 of the points
in $l\cap X$ by $x_1$, $x_2$ and $x_3$  and  $p=p_z(x_i)=l\cap \Bbb P^{n+1}$.
The hypersurface $p_z(X)\subset \Bbb P^{n+1}$ has at $p$ a  decomposition
into local irreducible components $p_z(X)\cap U_p=H_1\cup ... \cup H_k$,
where $U_p$ is a sufficiently small neighborhood of $p$. The points $x_i$ $i=1,...,3$
have neighborhoods $U_i$ such that $p_z:U_i \to p_z(U_i)$ is finite and
$p_z(U_i)$ contains one of $H_j$. Consider
the case where the local irreducible components $H_j$ contained by $p_z(U_i)$
are all distinct, w.l.o.g. denote them by $H_1$, $H_2$ and $H_3$ (the other cases
will follow by the same argument and are more favourable to our purposes). In this case
$H_1\cap H_2\cap H_3$ will be of dimension $n-2$. Since for every point
$t \in H_1\cap H_2\cap H_3$ the line $\overline {zt}$ is trisecant, the result
follows from corollary 1.11.

\

The part 2) follows from known facts about the trisecant varieties of
smooth varieties $X$ of codimension 2 in projective space $\Bbb P^N$.
The trisecant variety is irreducible if the dimension of $X$ $n\ge 2$.
Ziv Ran [Ra83] showed that if $n\ge 4$ and $Tr(X)\neq \Bbb P^{n+2}$ then
$X$ must be contained in a quadric (this result is not explicitly stated but
clearly follows from the article). Later Kwak [Ka02] showed that the same
holds for $n=3$. Ran also showed that if the degree of $X$ is less or equal to
its dimension then $X$, $d\le n$, then $X$ is a complete intersection. \par

We assume $X$ is nondegenerate in $\Bbb P^N$ (the degenerate case
follows the from hypersurface case). The above paragraph implies that if
$\dim Tr(X)=n+1$ then $C_XX=Tr(X)$ is the quadric containing $X$.  The
case $\dim Tr(X)=n$ or equivalently $Tr(X)=X$ is settled by a slicing argument
and the case $n=3$. It is  known, see for example remark 3.6 of [Ka02], that
if $X$ is of dimension $3$ and $Tr(X)=X$ then $X$ is a complete intersection of
two quadrics or the Segree variety $\Bbb P^1\times \Bbb P^2\subset \Bbb P^5$ which
is the intersection of 3 quadrics. If $n\ge 4$ consider a general $5$-plane
$L\subset \Bbb P^N$, then $X\cap L$ is a smooth 3-fold in $L=\Bbb P^5$ for which
$Tr(X\cap L)=X\cap L$, since $Tr(X\cap L)\subset Tr(X)\cap L$ and $X\cap L\subset
Tr(X\cap L)$. Then $X\cap L$ is one of the two cases described above. Both cases
have degree equal to 4 hence the degree $X$ is also 4. It follows from the result
of Ran and the end of the previous paragraph that $X$ is a complete intersection
of two quadrics.\hfill \ \qed
\enddemo

\

\

\proclaim {Theorem D} Let $X$ be a smooth subvariety of codimension 2
in $\Bbb P^{n+2}$. If $n\ge 3$ then:

$$H^0(X,S^m[\Omega^1_X(1)])=0$$

\noindent if and only if $X$ is not contained in a quadric.

\endproclaim

\demo {Proof} If $X$ is not contained in a quadric, then $C_XX=\Bbb P^{n+2}$
by proposition 1.12.  The vanishing $H^0(X,S^m[\Omega^1_X(1)])=0$ follows theorem C.\par
To analyse the case where $X$ is contained in a quadric  $Q$ recall that  proposition 1.5
states that $H^0(X, S^m[ {\Omega^1_X}(1)])=\{\Omega \in
S^m[\Bbb Cdz_0\oplus ...\oplus \Bbb Cdz_N]|
\text { } Z(\Omega|)\cap{T_xX})\text { is a cone with}$ $\text { vertex at } x,
\text { } \forall x \in X\}$. Consider the symmetric differential $\Omega_Q\in
S^2[\Bbb Cdz_0\oplus ...\oplus \Bbb Cdz_N]$ associated with the quadric $Q$.
For all $x\in X$ $Z(\Omega_Q)\cap{T_xX})(=Q\cap {T_xX})$  is a cone with vertex $x$ since $T_xX\in T_xQ$.
Hence $\Omega_Q$ defines an element of $H^0(X,S^m[\Omega^1_X(1)])$ and this element is
nontrivial since $Tan(X)=\Bbb P^{n+2}$.
\hfill \ \qed
\enddemo

\

\

\

\

\subhead  {1.3 Symmetric differentials on subvarieties of abelian
varieties }\endsubhead

\

\

In this section we do a short presentation of the results which are the
analogue to theorem A and part of theorem C for subvarieties of abelian
varieties. Again we are having in mind subvarieties with "low"
codimension. Recently, Debarre [De06] using the same perspective tackled
the problem of which subvarieties have an ample cotangent bundle, which
are in the other end in terms of codimension.

\

Let $X$ be a smooth subvariety of an abelian variety $A^n$. The surjection
on the conormal exact sequence:

$$0\to N_{X/A^{n}}^* \to \Omega^1_{A^n}|_{_X} \to \Omega^1_X \to 0 $$

\noindent induces the inclusion $j:\Bbb P(\Omega^1_X) \to
\Bbb P(\Omega^1_{A^n}|_{_X})$ of projectivized cotangent bundles.
The projectivized cotangent bundle of $A^n$ is trivial, i.e.
$\Bbb P(\Omega^1_{A^n})\simeq A^n\times\Bbb P^{n-1}$. Let
$p_2:\Bbb P(\Omega^1_{A^n})\to \Bbb P^{n-1}$ denote
the projection onto the second factor. Then $\Cal O_{\Bbb
P(\Omega^1_{A^n})}(m)\simeq p_2^*\Cal O_{\Bbb P^{n-1}}(m)$. The
composed map $f=p_2\circ j$:

$$f:\Bbb P(\Omega^1_{X})\to \Bbb P^{n-1}$$

\noindent is called the tangent map for $X$ in $A^n$.

\

\proclaim {Theorem F} Let $X$ be a smooth subvariety of an abelian variety $A^n$.
If the tangent map $f:\Bbb P(\Omega^1_{X})\to \Bbb P^{n-1}$ is both surjective and
connected then $\forall m\ge 0$:

$$H^0(X,S^m\Omega^1_X)=H^0(A^n,S^m\Omega^1_{A^n})$$.

\endproclaim \demo {Proof}  Associated to the tangent map $f$ is the  isomorphism:

$$\Cal O_{\Bbb
P(\Omega^1_{X})}(m)\simeq f^*\Cal O_{\Bbb P^{n-1}}(m) \tag 1.2.9$$

\noindent As in the proof of lemma 1, the isomorphism (1.2.9)
and the connectedness and surjectivity of the tangent map $f$ give that:

$$H^0(\Bbb P(\Omega^1_X),f^*(\Cal
O_{\Bbb P^{n-1}}(m)))=f^*H^0(\Bbb P^{n-1},\Cal O_{\Bbb P^{n-1}}(m))$$

\noindent The  result follows from the identifications $H^0(X,S^m\Omega^1_X)=
H^0(\Bbb P(\Omega^1_X),
\Cal O_{\Bbb P(\Omega^1_X) }(m))$ and $H^0(A^n,S^m\Omega^1_{A^n})=p_2^*
H^0(\Bbb P^{n-1},\Cal O_{\Bbb P^{n-1}}(m))$.
\hfill \ \qed
\enddemo

\

\proclaim {Corollary 1.13} Let $X$ be a smooth hypersurface of
an abelian variety $A^n$ with $n>2$  which does not contain any translate
of an abelian subvariety of $A^n$. Then $\forall m\ge 0$:

$$H^0(X,S^m\Omega^1_X)=H^0(A^n,S^m\Omega^1_{A^n})$$

\endproclaim \demo {Proof}  It follows from theorem F that it is
enough to show that the tangent map $f$ is connected and
surjective. The hypothesis on $X$ implies that $X$ itself is not the
translate of an abelian subvariety, see [Ab94] and hence the tangent map is
surjective.\par

The tangent map $f:\Bbb P(\Omega^1_{X})\to \Bbb P^{n-1}$ induces
the  map $\gamma:X \to \Bbb G(n-1,n)$, which is the Gauss map for
$X$ in $A^n$.  The fibers $f^{-1}(p)$
for $p\in \Bbb P^{n-1}$ project  to $F_p=\pi(f^{-1}(p))\subset X$.
The set $F_p$ consists of all the points  $x\in X$ for which the
line in $T_0A^n$ corresponding to $p$ is contained in $T_xX$. We
are using the common identification of the tangent space $T_xA^n$
for any $x\in A^n$ with $T_0A^n=\Bbb C^n$, which sends the tangent
spaces to $X$ at  $x\in X$  to an hyperplane of $T_0A^n$.
The Grassmanian $\Bbb G(n-1,n)$ is $\Bbb P^{n-1}$ and the
subvariety $W\subset \Bbb G(n-1,n)$ consisting of all $n-2$-planes
passing through the point $p \in \Bbb P^{n-1}$ is a hyperplane
$H \subset \Bbb G(n-1,n)$.\par

The hypothesis on $X$  guarantee  the the Gauss map
is finite (see corollary 3.10 of [Za93]). Hence the image of $X$
under the Gauss map is at least of dimension 2. It follows then that
$F_p=\gamma^{-1}(\gamma (X)\cap H)$ is connected by Bertini's theorem
and hence $f^{-1}(p)$ is also connected since $\pi:f^{-1}(p)\to F_p$
is 1 to 1.
\hfill \ \qed
\enddemo

\

As in the case of subvarieties of $\Bbb P^N$  we
obtain a vanishing theorem.

\

\proclaim {Theorem G} Let $X$ be a smooth subvariety of an abelian variety $A^n$.
If the general fiber of the tangent map $f:\Bbb P(\Omega^1_{X})\to \Bbb P^{n-1}$
is positive dimensional then $\forall m\ge 0$:

$$H^0(X,S^m\Omega^1_X\otimes L)=0$$

\noindent if $L$ is a negative line bundle on $X$.
\endproclaim

\demo {Proof} It follows from the arguments of theorem F and theorem A.
\hfill \ \qed
\enddemo

 \

 \

 \

\head {2. The non-invariance  of the cotangent plurigenera}\endhead

\

\

Let $X$ be a smooth projective variety. As in [Sa78], we define:

$$Q_m(X)=\dim H^0(X,S^m\Omega^1_X) \tag 2.1.1$$

\noindent The dimension $Q_m(X)$ is called the symmetric m-genus
of $X$. The graded ring $\Omega(X)=\sum_{m=0}^\infty
H^0(X,S^m\Omega^1_X)$ is called the cotangent ring of $X$. We define
the cotangent dimension of $X$ to be:

$$\lambda_I(X)=\dim_{Iitaka} \Cal O_{\Bbb P(\Omega^1_X)}(1) \tag 2.1.2$$

\noindent  where $\dim_{Iitaka} \Cal O_{\Bbb P(\Omega^1_X)}(1)$
is the Iitaka dimension of the the line bundle $\Cal O_{\Bbb P
(\Omega^1_X)}(1)$ on $\Bbb P(\Omega^1_X)$. For example, it follows
from the results of the previous section that if $X$ is a smooth
subvariety of $\Bbb P^N$ with $\dim_{\Bbb C}X>N/2$ then $\lambda_I(X)=
-\infty$. An abelian
variety $X$ of dimension $n$ has $\lambda_I(X)=n-1$ and a smooth
variety $Y$ of dimension $n$ with ample cotangent bundle has the maximal possible
Iitaka cotangent dimension for varieties of dimension $n$, $\lambda_I(Y)=2n-1$.

\

The symmetric 1-genus, $Q_1(X)$, of a smooth projective variety
$X$ is also called the irregularity of $X$. The irregularity of a
Kahler manifold is a topological invariant (as follows from Hodge
theory) and hence it can not jump in smooth families. One can also
see the irregularity of $X$  as one of the plurigenus of $X$, more
precisely $P_1(X)$ where $P_m(X)=\dim
H^0(X,(\bigwedge^n\Omega^1_X)^{\otimes m})$. There is an amazing
result of Siu [Si98] that states that all plurigenera are
invariant in smooth families of projective varieties.

\

The symmetric plurigenera behaves differently. The first author
gave the first example of a smooth family of projective varieties
where the cotangent m-genus jumps [Bo78]. We start by presenting
self contained modification of that example.

\

Let $T^3=\Bbb C^3/\Lambda$ be an abelian 3-fold, where $z_1$,
$z_2$ and $z_3$ are the Euclidean holomorphic coordinates of $\Bbb
C^3$. We denote the involution of $T^3$ given by the map $z\to
-z$ in $\Bbb C^3$ by $\sigma:T^3 \to T^3$. Let $X_t$ be a one-dimensional
family, over the disc $\Delta$, of $\sigma$-invariant smooth hypersurfaces of
$T^3$ which pass through one of the fixed points $p_0$ of the involution
$\sigma$. Moreover, locally on a neighborhood of the fixed
point $p_o$, which we assume to be $p_0=(0,0,0)$,  $X_t$ is given by the
equation:

$$z_1=tz^3_2+f_t(z_2,z_3) \tag 2.1.3$$

\noindent where $f_t(z_2,z_3)\in (z_2,z_3)^5$.

\

\proclaim {Theorem H}  Let $X_t$ be the family described in (2.1.3) and
$Y_t$ be the family which is the simultaneous minimal resolution
of the  family nodal varieties $V_t=X_t/\sigma$. Then: \smallskip  a) The
symmetric plurigenera is not invariant along the family
$Y_t$.\smallskip b) Moreover $H^0(Y_0,S^m\Omega^1_{Y_0})
\simeq[\sum_{m_1\ge m_2+m_3} \Bbb C dz_1^{m_1}dz_2^{m_2}dz_3^{m_3}]^{\Bbb Z_2}$ and
$H^0(Y_t,S^m\Omega^1_{Y_t})\supset [\sum_{3m_1\ge m_2+m_3}
\Bbb C dz_1^{m_1}dz_2^{m_2}dz_3^{m_3}]^{\Bbb Z_2}$
\endproclaim

\demo {Proof} We have, first, to describe the local picture for
the symmetric differentials around the fixed point $p_0$. Let
$(U,0)$ be the neighborhood germ of the origin in $\Bbb C^2$ and
$(U/\sigma,x_0)$ be the neighborhood germ of the  nodal surface
singularity ($\sigma$ is again the -id involution).  Let $(V,E)$ be the
neighborhood germ of the (-2)-curve and $r:(V,E) \to
(U/\sigma,x_0)$ the minimal resolution.

\

\proclaim {Lemma 2.1} The symmetric differentials on $(V,E)$ are in
bijection with the symmetric differentials on $(U,0)$ which are
$\sigma$-invariant and of the form $\omega=\sum
h_{m_1,m_2}(z_1,z_2)dz_1^{m_1}dz_2^{m_2}$, where
$h_{m_1,m_2}(z_1,z_2)\in (z_1,z_2)^{m_1+m_2}\subset \Cal O(U)$.
\endproclaim

\demo {Proof} Let $(W,E')$ be the germ neighborhood of the blow up
of $(U,0)$ at $0$ and $b:(W,E')\to (U,0)$ the blow up map. Denote
by $g:(W,E')\to (V,E)$ be the 2 to 1 naturally defined covering of
$V$ ramified at $E\subset V$ (for which  $\sigma\circ b=r\circ g$
holds). \par First, we note that there is a natural bijection
between $H^0(V\setminus E,S^m\Omega^1_V)$  with  $[H^0(W\setminus
E',S^m\Omega^1_{W})]^{\Bbb Z_2}$. The differential pullback $dg^*:
H^0(V,S^m\Omega_V^1) \to H^0(W,S^m\Omega^1_W)$ is an injection.
We want to see the differentials of $W$ on $U$, to do this we notice
that there is a natural identification $H^0(U,S^m\Omega^1_U)=
H^0(W,S^m\Omega^1_W)$. Hence, we have that
$H^0(V,S^m\Omega_V^1)\subset [H^0(U,S^m\Omega^1_{U})]^{\Bbb Z_2}$.
What remains to be determined  is  what to require on the
elements of $[H^0(U,S^m\Omega^1_U)]^{\Bbb Z_2}$ for them to be, after
the natural identification mentioned above, pullbacks by $dg^*$ .

We give a coordinate chart approach to this problem. Consider the
single affine chart blow up map $f:(\Bbb C^2,u,v) \to (\Bbb C^2,z_1,z_2)$ given by
$(u,v) \to (u,uv)$. Consider also the double ramified covering
$g:(\Bbb C^2,u,v) \to (\Bbb C^2,x,y)$ given by $(u,v)\to
(u^2,v)$.

\

The relations between the differentials in the different
coordinate charts $\Bbb C^2$ are:

$$f^*(dz_1)=du, \text { }\text { }
f^*(dz_2)=udv+vdu \text { }  \text { }  \text { } \text { }
  \text {;}
\text { }
  \text { }    \text { }
  \text { }  \text { }     g^*(dx)=2udu,  \text { } \text { }
g^*(dy)=dv \tag 2.13$$

Let us write the pullback by $f$ of a symmetric differential on
$(\Bbb C^2,z_1,z_2)$:

$$f^*(dz_1^{m_1}dz_2^{m_2})=\sum_{i=0}^{m_2} {m_2\choose i}
u^{m_2-i}v^idv^{m_2-i}du^{m_1+i} \tag 2.1.4$$

\noindent  From (2.1.4) one easily sees that the pullback by $f$ of any
symmetric differential monomial of order $m$ ($m_1+m_2=m$) has the term
in $du^{m}$ with no power of $u$ in the coefficient (and all other terms have
$du$ with an order smaller than $m$). On the other hand, by (2.1.3)
a symmetric differential $\omega=u^{i_1}v^{i_2}du^{m_1}dv^{m_2}$
is a pullback by $dg^*$ only if $i_1\ge m_1$. This  plus the coordinate
description of $f$ implies that
 $f^*(h(z_1,z_2)dz_1^{m_1}dz_2^{m_2})$ with $m_1+m_2=m$
is  a pullback by $g$ only if the Taylor expansion of $h$ has
all terms with combined order in $z_1$ and $z_2$ to be greater or
equal to $m$, which concludes the proof.

\hfill \ \qed

\enddemo

Let us set up some notation. The maps $q_t:X_t \to X_t/\sigma$ and
$r_t:Y_t \to X_t/\sigma$ will denote respectively the quotient map
induced by $\sigma$ and the minimal resolution of $V_t$. Consider
the auxiliary family $W_t$ whose members are the $X_t$ with $p_0$
blown up, denote the blowing up map by $b_t:W_t \to X_t$. Applying
the same argument of the lemma 2.1 it follows that there is a
naturally defined injection of $j:H^0(Y_t,S^m\Omega^1_{Y_t})
\hookrightarrow H^0(X_t,S^m\Omega^1_{X_t})$.
\par It follows from corollary 1.10  that the m-the order
symmetric differentials  on $X_t$ are the restrictions to $X_t$ of
the symmetric differentials $\omega$ of the abelian variety $T^3$
($T^3$ can be chosen to be simple), which are the form:

$$\omega=\sum_{_{m_1+m_2+m_3=m}}
a_{_{m_1,m_2,m_3}}dz_1^{m_1}dz_2^{m_2}dz_3^{m_3}$$

\noindent where $a_{_{m_1,m_2,m_3}}\in \Bbb C$. We want to describe
the symmetric differentials $j(H^0(Y_t,S^m\Omega^1_{Y_t}))
\subset H^0(X_t,S^m\Omega^1_{X_t})$. Recall that the hypersurfaces $X_t$
 are given in a sufficiently small neighborhood of $p_0$ by
the equation $z_1=tz^3_2+f_t(z_2,z_3)$, where $f_t(z_2,z_3)\in
(z_2,z_3)^5$  and $z_1$,
$z_2$ and $z_3$ are the Euclidean holomorphic coordinates of $\Bbb
C^3$. We can use $z_2$ and $z_3$ as local holomorphic
coordinates of a sufficiently small neighborhood $U_t\subset X_t$
of $p_0$. Let $i_t:X_t \hookrightarrow T^3$ be the inclusion map,
then the pullback  of $dz_1$ to
$U_t\subset X_t$ is:

$$i_t^*(dz_1)|_{U_t}=3tz_2^2dz_2+\mu_t \tag 2.1.5$$

where $\mu_t\in (z_2,z_3)^4dz_2+(z_2,z_3)^4dz_3$ and hence:

$$i_t^*(dz_1^{m_1}dz_2^{m_2}dz_3^{m_3})|_{U_t}=
ct^{m_1}z_2^{2m_1}dz_2^{m_1+m_2}dz_3^{m_3}+\gamma_t \tag 2.1.6$$

\noindent where the $(m_1+m_2+m_3)$-order symmetric differential $\gamma_t$
has the coefficients of their monomial terms in $(z_2,z_3)^{2m_1+2}$.
It follows from the Lemma 2.1 that if
$t\neq 0$ the symmetric differential
$i_t^*(dz_1^{m_1}dz_2^{m_2}dz_3^{m_3}) \in
j(H^0(Y_t,S^m\Omega^1_{Y_t}))$ if and only if $2m_1\ge m_1+m_2+m_3$ and
$m_1+m_2+m_3=0$ mod 2, giving the first part of b).

 \

  If $t=0$ we have
$i_0^*(dz_1^{m_1}dz_2^{m_2}dz_3^{m_3})|_{U_t}=\mu_t^{m_1}dz_2^{m_2}dz_3^{m_3}$.
The $(m_1+m_2+m_3)$-order symmetric differential $\mu_t^{m_1}dz_2^{m_2}dz_3^{m_3}$
has the coefficients of their monomial terms in $(z_2,z_3)^{4m_1}$.
As before, it follows from the lemma 2.1 that the
symmetric differential $i_0^*(dz_1^{m_1}dz_2^{m_2}dz_3^{m_3}) \in
j(H^0(Y_0,S^m\Omega^1_{Y_0}))$ if and only if $4m_1\ge m_1+m_2+m_3$ and
$m_1+m_2+m_3=0$ mod 2, which gives the second part of b) and concludes the proof.

\hfill \ \qed\enddemo

\

\

We now answer a question posed by Paun: are the dimensions of
$H^0(X_t,S^m\Omega^1_{X_t}\otimes  K_{X_t})$ constant for a family
of smooth projective varieties? We answer this question
negatively.

\

\proclaim {Theorem I} Let $Y_t$ be a family of smooth projective
varieties. The invariance of the dimension of
$H^0(Y_t,S^m\Omega^1_{Y_t}\otimes K_{Y_t})$ does not necessarily
hold along the family.
\endproclaim

\demo {Proof} Let $X_t$ be a family over $\Delta$ of smooth
hypersurfaces of degree $d$ of $\Bbb P^3$ specializing to a nodal
hypersurface $X_0$ with $l>\frac {8}{3}(d^2-\frac{5}{2}d)$ nodes.
This is possible as long $d\ge 6$ [Mi83]. Denote by $Y_t$ the
family  which is the simultaneous resolution of the family $X_t$,
$t\in \Delta$. The general member of $Y_t$ is a smooth
hypersurface of $\Bbb P^3$ of degree $d$ and $Y_0$ is the minimal
resolution of $X_0$. We proved in [BoDeO06] that if $d\ge 6$ $Y_0$
has plenty of symmetric differentials, more precisely
$H^0(Y_0,S^m\Omega^1_{Y_0})\uparrow m^3$. This result plus the
effectivity of the canonical divisor $K_{Y_0}$, in particular,
implies that:

$$ H^0(Y_0,S^m\Omega^1_{Y_0}\otimes K_{Y_0})\neq 0, \text { } \text {
} m\gg 0 \tag 2.1.7$$
\par Theorem B  gives that $H^0(Y_t,S^m\Omega^1_{Y_t}\otimes
\Cal O_{Y_t}(m))=0$. The canonical divisor of the hypersurface
$Y_t$ is $K_{Y_t}=\Cal O_{Y_t}(d-4)$, which implies that :

$$H^0(Y_t,S^m\Omega^1_{Y_t}\otimes K_{Y_t})=0, \text { } \text { } m\ge
d-4 \tag 2.1.8$$

\noindent The result follows from (1.2.7) and (1.2.8). \hfill \
\qed\enddemo

\

We want to modify the question of Paun to be able to express a
stronger result (which is suggested by the proof of the above
theorem). We introduce the following notation:

$$Q_{\alpha,m}(X)=\dim H^0(X,S^m(\Omega^1_X\otimes \alpha K_X)) \tag 2.1.9$$

The dimension $Q_{\alpha,m}(X)$ is called the $\alpha$-twisted
symmetric m-genus of $X$. Is there an $\alpha$ for which the
$\alpha$-twisted plurigenera is invariant along all families of
smooth projective varieties $X_t$ with $K_{X_t}>0$?

\

\proclaim {Question} What is the lower bound $\beta$ for the
$\alpha$'s for which the $\alpha$-twisted plurigenera is invariant
along all families of smooth projective varieties $X_t$ with $K_{X_t}>0$?
\endproclaim

\

\proclaim {Proposition 2.2} The lower bound $\beta$ for the $\alpha$'s
for which the $\alpha$-twisted plurigenera is invariant for the
families of smooth projective surfaces $X_t$ with $K_{X_t}>0$ must satisfy:

$$\beta\ge 1/2$$
\endproclaim

\demo {Proof} Consider the families $Y_t$ described in theorem 1
for degree $d=6$. In this case $K_{Y_t}=2\Cal O_{Y_t}(1)$ and by
theorem B we have $H^0(X_t,S^m(\Omega^1_{Y_t}\otimes 1/2
K_{Y_t}))=0$ for $t\neq 0$. On the other hand,
$H^0(Y_0,S^m(\Omega^1_{Y_0}\otimes 1/2 K_{Y_0}))\neq 0$ for $m\ge
2$ as in  theorem I. \hfill \ \qed\enddemo

\

\

\end
\Refs

\ref \key Ab94  \by {\text { }\text { }\text { }\text { }\text {
}\text { }\text { }\text { }\text { }\text { }\bf D.Abramovich} \paper
 Subvarieties of semiabelian varieties \jour Comp. Math.   \vol 90 (1)
\yr 1994 \pages 37--52
\endref

\ref  \key Bo1-78 \by {\text { }\text { }\text { }\text { }  \text {
}\text { }\text { }\text { }\text { }\bf F.Bogomolov} \paper
Unstable vector bundles and curves on surfaces\jour Proc. of the
Int. Cong. of \text { }\text { }\text { }\text {
}\text { }\text { }\text { }\text { }\text { }\text { }Mathematicians
(Helsinki, 1978) \pages 517-524
\endref

\ref  \key Bo2-78 \by {\text { }\text { }\text { }\text { }  \text {
}\text { }\text { }\text { }\text { }\bf F.Bogomolov} \paper
Holomorphic symmetric tensors on projective surfaces \jour Uspekhi
Mat. Nauk \text { }\text { }\text { }\text { }\text { }\text {
}\text { }\text { }\text { }\vol 33 (5) \yr 1978 \pages 171 - 172
 (Russian Math. Surveys v 33 (5) (1978) p 179 -180)
\endref

\ref  \key BoDeO05 \by {\text { }\text { }\text { }\text { } \text
{ }\text { }\text { }\text { }\text { }\bf F.Bogomolov, B.De
Oliveira} \paper Hyperbolicity of nodal hypersurfaces \jour J.
Regne Angew. \text { }\text { }\text { }\text { }\text { }\text {
}\text { }\text { }\text { }\text { }Math.  \vol 596 \yr 2006\pages
\endref

\ref  \key Br71 \by {\text { }\text { }\text { }\text { }  \text {
}\text { }\text { }\text { }\text { }\bf P.Bruckmann} \paper
Tensor differential forms on algebraic varieties \jour Izvestya
A.N. USSR \vol 35 \text { }\text { }\text { }\text { }\text {
}\text { }\text { }\text { }\yr 1971 \pages (English translation
Math USSR Izvestya 1979 v 13 (1) p 499 -544)
\endref


\ref \key De05  \by {\text { }\text { }\text { }\text { }\text {
}\text { }\text { }\text { }\text { }\text { }\bf O.Debarre} \paper
Varieties with ample cotangent bundle \jour Comp. Math.   \vol 141
\yr 2005 \pages 1445-1459
\endref

\ref \key Ka02  \by {\text { }\text { }\text { }\text { }\text {
}\text { }\text { }\text { }\text { }\text { }\bf S. Kwak} \paper
Smooth threefolds in $\Bbb P^5$ without apparent triple or quadruple
points and a \text { }\text { }\text { }\text {
}\text { }\text { }\text { }\text { }\text { }quadruple-point formula  \jour Math. Ann. \vol 230 \yr 1983
\pages 649-664
\endref

\ref \key Ra83  \by {\text { }\text { }\text { }\text { }\text {
}\text { }\text { }\text { }\text { }\text { }\bf Z. Ran} \paper
On projective varieties of codimension 2  \jour Inv. Math. \vol 73 \yr 1983
\pages 333-336
\endref

\ref \key RaLo03\by {\text { }\text { }\text { }\text { }\text {
}\text { }\text { }\text { }\text { }\text { }\bf Z. Ran; A.Lopez}\paper
On the irreducibility of secant cones, and an application to linear\text { }
\text { }\text { }\text {
}\text { }\text { }\text { }\text { }\text { } normality
\jour Duke Math. J.\vol 117 (3) \yr 2003
\pages 389-401
\endref


\ref \key Sa79  \by {\text { }\text { }\text { }\text { }\text {
}\text { }\text { }\text { }\text { }\text { }\bf F.Sakai} \paper
Symmetric powers of the cotangent bundle and classification of
algebraic varieties \jour \text { }\text { }\text { }\text {
}\text { }\text { }\text { }\text { }\text { } Lect. Notes in
Math. (Springer Verlag) \vol 732 \yr 1979 \pages 545-563
\endref

\ref \key Sc92  \by {\text { }\text { }\text { }\text { }\text {
}\text { }\text { }\text { }\text { }\text { }\bf M.Schneider} \paper
Symmetric differential forms as embedding obstructions and vanishing\text { }
\text { }\text { }\text { }\text { }\text { }\text { }
\text { } \text { }\text { }\text { }\text { }\text { }\text { }\text {
}\text { }theorems  \jour J.Algebraic Geom. \vol 1 \yr 1992
\pages 175-181
\endref


\ref \key Si98  \by {\text { }\text { }\text { }\text { }\text {
}\text { }\text { }\text { }\text { }\text { }\bf Y-T.Siu} \paper
Invariance of plurigenera  \jour Inv. Math. \vol 134 \yr 1998
\pages 661-673
\endref

\ref  \key Za81 \by {\text { }\text { }\text { }\text { }  \text {
}\text { }\text { }\text { }\text { }\bf F.Zak} \paper Projections
of algebraic varieties \jour Izvestya Mat. Sb. \vol 116 \text {
(3)} \yr 1981 \pages 593-602 \text { }\text { }\text { }\text { }
\text { } \text { }\text { }\text { }\text { }\text { }\text {
}\text { }  \text { } (English translation Math USSR Izvestya 1979
v 13 (1) p 499 -544)
\endref

\ref  \key Za93 \by {\text { }\text { }\text { }\text { }  \text {
}\text { }\text { }\text { }\text { }\bf F.Zak} \paper Tangents and
secants of algebraic varieties \jour Trans. Math. Mon. AMS \vol 127
 \yr 1993
\endref
\endRefs
